\documentclass[12pt]{article}
\usepackage{amsmath,amssymb}

\newtheorem{de}{Definition}[section]
\newtheorem{lm}[de]{Lemma}
\newtheorem{pr}[de]{Proposition}

\newtheorem{re}[de]{Remark}

\newtheorem{te}[de]{Theorem}

\begin{document}

\title{\textbf{From regular modules to von Neumann regular rings via coordinatization}\thanks{Research supported by the National Research Foundation, United Arab Emirates University (NRF-UAEU) under Grant No. 31S07.}}

\author{\textbf{Leonard D\u{a}u\c{s}$^{1,2}$ and Mohamed A. Salim$^{3}$} \\
{\small $^{1,3}$Department of Mathematical Sciences, United Arab Emirates University,}\\
{\small Al Ain, P.O.Box 15551, UAE}\\
{\small $^{2}$Departament of Mathematics, Technical University of Civil Engineering,}\\
{\small Bdul. Lacul Tei 124, RO-020396 Bucharest 2, Romania}\\
{\small $^{1}$leonard.daus@uaeu.ac.ae, $^{2}$daus@utcb.ro, $^{3}$msalim@uaeu.ac.ae}}


\renewcommand{\thefootnote}{\fnsymbol{footnote}}    
\footnotetext{\emph{2010 Mathematics Subject Classification:} 16D80, 16E50, 06C20, 06D99} 
\footnotetext{\emph{Keywords:} coordinatization, regular module, von Neumann regular ring, complemented\\ modular lattice, distributive lattice} 
\renewcommand{\thefootnote}{\arabic{footnote}}
\date{}
\maketitle

\begin{abstract}
In this paper we establish a very close link (in terms of von Neumann's coordinatization) between regular modules introduced by Zelmanowitz, on one hand, and von Neumann regular rings, on the other hand: we prove that the lattice $\mathcal{L}^{fg}(M)$ of all finitely generated submodules of a finitely generated regular module $M$, over an arbitrary ring, can be coordinatized as the lattice of all principal right ideals of some von Neumann regular ring $S$. 
\end{abstract}


\section{Introduction and Main Results}

${\;\;\;\;\;\;\;}$In 1936, John von Neumann defined a ring $R$ to be regular if for any $r\in R$ there exists $s\in R$ such that $r=rsr$. Motivated by the coordinatization of projective geometry, which was being reworked at that time in terms of lattice, von Neumann introduced regular rings as an algebraic tool for studying certain lattice: in fact, regular rings were used to coordinatize complemented modular lattice, a lattice $L$ being coordinatized by a regular ring $R$ if it is isomorphic to the lattice $\mathbf{L}(R_{R})$ of all principal right ideals of $R$.\newline
${\;\;\;\;\;\;\;}$J. Zelmanowitz in paper \cite{z} followed the original elementwise definition of von Neumann and called a right $R-$module $M$ regular if for any $m\in M$ there exists $g\in Hom_{R}(M,R)$ such that $mg(m)=m$. Since a morphism $f\in Hom_{R}(R,M)$ is uniquely given by an element $m\in M$, one can reformulate the regular module defined by Zelmanowitz as follows: for any $f\in Hom_{R}(R,M)$ there exists $g\in Hom_{R}(M,R)$ such that $f=f\circ g\circ f$. \newline
${\;\;\;\;\;\;\;}$In paper \cite{dntd} it was defined the concept of a regular object with respect to another object (or a relative regular object) in an arbitrary category, which extends the notion of regular module. In particular, when we consider the category $\mathcal{M}_{R}$ of right $R-$modules we obtain the following definition
\begin{de}
Let $M$ and $U$ be two right $R-$modules. We say that $M$ is $U-$regular module if for any $f\in Hom_{R}(U,M)$ there exists a
morphism $g\in Hom_{R}(M,U)$ such that $f=f\circ g\circ f$.
\end{de}
${\;\;\;\;\;\;\;}$Obviously, $M$ is a regular module if and only if $M$ is $R-$regular. The concept of relative regular module has been proved to be an extremely useful tool in the theory of von Neumann regular rings and regular modules (see \cite{dntd} and \cite{d}).\newline
${\;\;\;\;\;\;\;}$For a right $R-$module $M$ we denote by $\mathcal{L}^{fg}(M)$ the lattice of all finitely generated submodules of $M$, partially ordered by inclusion. The aim of this paper is to study the lattice $\mathcal{L}^{fg}(M)$ of all finitely generated submodules of a finitely generated regular module. We first prove that $\mathcal{L}^{fg}(M)$ is a complemented modular lattice. 

\begin{de}
A lattice $\mathcal{L}$ is coordinatizable if there exists a von Neumann regular ring $R$ such that $\mathcal{L}$ is isomorphic to the lattice $\mathbf{L}(R_{R})$ of all principal right ideals of $R$.
\end{de}

${\;\;\;\;\;\;\;}$Using the isomorphism between the lattices $\mathcal{L}^{fg}(M)$ and $\mathcal{L}^{fg}(S_{S})$, where $S$ is the endomorphisms ring $End_{R}(M)$, we obtain our main result:

\begin{te}
If $R$ is an arbitrary ring and $M$ is a finitely generated regular right $R-$module, then $\mathcal{L}^{fg}(M)$ is coordinatizable.
\end{te}

${\;\;\;\;\;\;\;}$The key step in proving this theorem is \emph{Proposition 2.3}, which was inspired by a result given by von Neumann, for the case $R_{R}$, in his classical paper \cite{vn}. In order to adapt his result in the general case of regular modules we use the concept of relative regular module. \newline
${\;\;\;\;\;\;\;}$A ring $R$ is called strongly regular if for any element $r\in R$ there exists an element $s\in R$ such that $r^2s=r.$ Finally, we will prove:

\begin{te}
Let $M$ be a finitely generated regular right module over an arbitrary ring $R$. Then $\mathcal{L}^{fg}(M)$ is distributive if and only if $End_{R}(M)$ is strongly regular ring.
\end{te}
 
${\;\;\;\;\;\;\;}$In this paper, all modules will be right modules. We refer the reder to Birkhoff \cite{b} or Gratzer \cite{gr} for elements of lattice theory and to \cite{g} for definitions and results about von Neumann regular rings.

\section{Proofs of the Main Results}

In order to prove $\mathcal{L}^{fg}(M)$ is a complemented modular lattice, we will need the following preliminary result:

\begin{lm}
Let $M$ be a regular $R-$module, let $U$ be a finitely generated $R-$module and let $f\in Hom_{R}(U,M)$. If $N$ is a finitely generated submodule of $M$, then $f^{-1}(N)$ is finitely generated $R-$module.
\end{lm}

\emph{Proof.}
Since $M$ is a regular module and $N$ is a finitely generated submodule of $M$, by \cite[Theorem 1.6]{z} we obtain that $N$ is a direct summand of $M$. Therefore there exists a submodule $N'$ of $M$ such that $M=N\oplus N'$. But $M/N=M+N'/N \simeq N'/N\cap N'\simeq N'$ and then $M/N$ is isomorphic with a direct summand of $M$, so $M/N$ is a regular module. Since $M/N$ is a regular module and $U$ is a finitely generated module, by \cite[Remark 4.3]{dntd} it follows that $M/N$ is $U-$regular.\newline
${\;\;\;\;\;\;\;}$Let $\pi:M\longrightarrow M/N$ be the canonical projection and we consider the morphism of $R-$modules $\pi\circ f:U\longrightarrow M/N$. By \cite[Proposition 3.1]{dntd}, we obtain that $Ker(\pi\circ f)=f^{-1}(N)$ is a direct summand of $U$. Obviously, $f^{-1}(N)$ is a finitely generated $R-$module.
\newline

\begin{pr}
Let $R$ be an arbitrary ring and let $M$ be a finitely generated regular right $R-$module. Then $\mathcal{L}^{fg}(M)$ is a complemented modular lattice $($where $N_1\vee N_2=N_1+N_2$ and $N_1\wedge N_2=N_1\cap N_2$, for any $N_1,N_2 \in \mathcal{L}^{fg}(M))$.
\end{pr}

\emph{Proof.}
Let $N_1,N_2 \in \mathcal{L}^{fg}(M)$. Obviously, $N_{1}+N_{2}\in \mathcal{L}^{fg}(M)$ and $N_1\vee N_2=N_1+N_2$. We want to prove that $N_1\cap N_2\in \mathcal{L}^{fg}(M)$. We consider the embedding morphism $i:N_{1}\hookrightarrow M$. Then $N_1\cap N_2=i^{-1}(N_{2})$, which is finitely generated by Lemma 2.1. Since $N_1\wedge N_2=N_1\cap N_2$, we may assert that $\mathcal{L}^{fg}(M)$ is a lattice.\newline
${\;\;\;\;\;\;\;}$In general, the lattice $\mathcal{L}(M)$ of all submodules of an arbitrary module $M$ with the lattice operations $+$ and $\cap$ is modular. Thus, $\mathcal{L}^{fg}(M)$ is a modular lattice, when $M$ is a regular module.\newline
${\;\;\;\;\;\;\;}$If $N\in \mathcal{L}^{fg}(M)$, by \cite[Theorem 1.6]{z} we obtain that $N$ is a direct summand of $M$. Hence we may find a submodule $N'$ of $M$ such that $M=N\oplus N'$. Since $N'$ is a direct summand of $M$ and $M$ is finitely generated as $R-$module, then $N'$ is finitely generated and therefore $N'$ is a complement of $N$ in $\mathcal{L}^{fg}(M)$. Thus $\mathcal{L}^{fg}(M)$ is a complemented lattice. 
\newline

\begin{pr}
Let $M$ be a finitely generated regular right $R-$module and $S=End_{R}(M)$. Then the lattices $\mathcal{L}^{fg}(M)$ and $\mathcal{L}^{fg}(S_{S})$ are isomorphic.
\end{pr}

\emph{Proof.}
Since $M$ is a finitely generated regular $R-$module, by \cite[Corollary 4.2]{z} we obtain that $S=End_{R}(M)$ is a von Neumann regular ring. We consider an arbitrary $J\in \mathcal{L}^{fg}(S_{S})$. Then there exists an idempotent $e\in S$ such that $J=eS$.It follows that $JM=eSM=eM\in \mathcal{L}^{fg}(M)$. Thus we obtain the monotone map
\[
\varphi:\mathcal{L}^{fg}(S_{S})\longrightarrow \mathcal{L}^{fg}(M), \ \ \ \varphi(J)=JM.
\]
We will prove that $\varphi$ is a lattice isomorphism.\newline
${\;\;\;\;\;\;\;}$For $N\in \mathcal{L}^{fg}(M)$ we denote $\psi(N)=\left\{f\in S/ f(M)\leq N\right\}$. Since $N$ is a direct summand of $M$, then there exists an idempotent $e\in S$ such that $e(M)=N$. Of course, $e_{|_N}=id_{N}$. We show that $\psi(N)=eS$.\newline
${\;\;\;\;\;\;\;}$Consider $f\in S$. It follows that $f(M)\leq M$ and thus $e(f(M))\leq e(M)=N$. Hence $ef\in \psi(N)$ and therefore $eT\subseteq \psi(N)$. For the converse inclusion we consider $g\in \psi(N)$. Then $g\in S$ and $g(M)\leq N$. Since $e(g(m))=g(m)$, for all $m\in M$, we obtain that $g\in eS$, so $\psi(N)\subseteq eS$. Hence $\psi(N)=eS\in \mathcal{L}^{fg}(S_{S})$ and therefore we obtain the monotone map
\[
\psi:\mathcal{L}^{fg}(M)\longrightarrow \mathcal{L}^{fg}(S_{S}), \ \ \ \psi(N)=\left\{f\in S/ f(M)\leq N\right\}
\]
${\;\;\;\;\;\;\;}$On one hand, if $e\in S$ is an arbitrary idempotent, then $\varphi(eS)=eM$ and $\psi(\varphi(eS))=\psi(eM)=eS$, so $\psi\circ\varphi=id_{\mathcal{L}^{fg}(S_{S})}$. On the other hand, $\varphi(\psi(eM))=\varphi(eS)=eM$, so $\varphi\circ\psi=id_{\mathcal{L}^{fg}(M)}$. Thus, we can conclude that $\varphi$ and $\psi$ are lattice isomorphisms and therefore the lattices $\mathcal{L}^{fg}(M)$ and $\mathcal{L}^{fg}(S_{S})$ are isomorphic.
\newline

\begin{re}
If $R$ is a von Neumann regular ring, by \cite[Theorem 1.1]{g} we obtain that the lattice $\mathcal{L}^{fg}(R_{R})$ is, in fact, the lattice $\mathbf{L}(R_{R})$ of all principal right ideals of $R$. 
\end{re}

We are now in the position to prove our main results:\newline

\emph{Proof of Theorem 1.3.}
When $M$ is a finitely generated regular module, over an arbitrary ring $R$, by \cite[Corollary 4.2]{z} it follows that $S=End_{R}(M)$ is a von Neumann regular ring. Our Theorem is now a direct consequence of \emph{Definition 1.2}, \emph{Proposition 2.3} and \emph{Remark 2.4}.
\newline

\emph{Proof of Theorem 1.4.}
Utumi proved in \cite[Theorem 1.1]{u} the following: \newline
A regular ring $R$ is strongly regular if and only if the lattice $\mathbf{L}(R_{R})$ of all principal right ideals of $R$ is distributive.\newline
${\;\;\;\;\;\;\;}$Using the lattice isomorphism 
\[
\mathcal{L}^{fg}(M)\simeq \mathcal{L}^{fg}(S_{S})= \mathbf{L}(S_{S})
\]
our theorem is now obvious.

\end{document}